\pdfoutput=1
\documentclass[letterpaper]{amsart}
\usepackage{amsthm}
\usepackage{amsmath,amssymb}
\usepackage{amsaddr}
\usepackage{xurl}
\usepackage{hyperref}
\usepackage{cleveref}
\hypersetup{colorlinks=true,citecolor=red}

\newtheorem{Th}{Theorem}
\newtheorem{Pro}[Th]{Proposition}
\newtheorem{Lem}[Th]{Lemma}
\newtheorem{Cor}[Th]{Corollary}

\theoremstyle{definition}
\newtheorem{Def}[Th]{Definition}
\newtheorem{Exp}{Example}
\newtheorem{Rem}{Remark}
\newtheorem{Prob}{Problem}

\theoremstyle{definition}

\newenvironment{Proof}{\begin{proof}}{\end{proof}}
\newcommand{\prover}{Prover}
\newcommand{\delayer}{Delayer}
\newcommand\gen[1]{\langle{#1}\rangle}
\DeclareMathOperator{\Nil}{Nil}
\DeclareMathOperator{\Jac}{Jac}
\DeclareMathOperator{\Reg}{Reg}
\newcommand\N{\mathbb{N}}
\newcommand\Z{\mathbb{Z}}
\newcommand\R{\mathbb{R}}

\title{A quantitative general Nullstellensatz for Jacobson rings}
\author{Ryota Kuroki}

\email{kuroki-ryota128@g.ecc.u-tokyo.ac.jp}
\address{Graduate School of Mathematical Sciences, The University of Tokyo, 3-8-1 Komaba, Meguro-ku, Tokyo, 153-8914, Japan}
\date{}

\begin{document}

\begin{abstract}
    The general Nullstellensatz states that if $A$ is a Jacobson ring, $A[X]$ is Jacobson.
    We introduce the notion of an $\alpha$-Jacobson ring for an ordinal $\alpha$ and prove a quantitative version of the general Nullstellensatz: if $A$ is an $\alpha$-Jacobson ring, $A[X]$ is $(\alpha+1)$-Jacobson.
    The quantitative general Nullstellensatz implies that $K[X_1,\ldots,X_n]$ is not only Jacobson but also $(1+n)$-Jacobson for any field $K$.
    It also implies that $\mathbb{Z}[X_1,\ldots,X_n]$ is $(2+n)$-Jacobson.
\end{abstract}

\maketitle

\section{Introduction}
In this paper, all rings are considered to be commutative with identity.
This paper contains non-constructive arguments.

In classical mathematics, a ring $A$ is called Jacobson if every prime ideal of $A$ is an intersection of maximal ideals.
One of the most important theorems about Jacobson rings is the general Nullstellensatz: a univariate polynomial ring over a Jacobson ring is Jacobson. This theorem has been independently proved by Goldman \cite[Theorem 3]{Go51} and Krull \cite[Satz 1]{Kr51}.

In \cite[Theorem 3.9]{Ku24}, we have extracted a constructive proof of the general Nullstellensatz from a classical proof \cite[Theorem 8]{Em} and revealed the computational content of the theorem.
In this paper, based on the computational content that we have obtained, we introduce the notion of an $\alpha$-Jacobson ring for an ordinal $\alpha$ and prove a quantitative version of the general Nullstellensatz: if $A$ is an $\alpha$-Jacobson ring, $A[X]$ is $(\alpha+1)$-Jacobson (corollary \ref{qgn}).

Informally speaking, a ring $A$ is called $\alpha$-Jacobson if we can prove that $A$ is Jacobson by an \emph{elementary proof of complexity at most $\alpha$}.
If a ring $A$ is $\alpha$-Jacobson for some $\alpha$, then $A$ is Jacobson (corollary \ref{a-jac-jac}). The converse does not hold (remark \ref{not-a-jac}).
Every $0$-dimensional ring is $1$-Jacobson (proposition \ref{0-dim-jac}), and $\Z$ is $2$-Jacobson (proposition \ref{1-dim-jac}).

Although the constructive general Nullstellensatz \cite[Theorem 3.9]{Ku24} does not follow from the quantitative general Nullstellensatz (corollary \ref{qgn}), the theorem provides a stronger result for typical examples. For example, we can prove that $\Z[X_1,\ldots,X_n]$ is not only Jacobson but also $(2+n)$-Jacobson.

In section 4, we introduce the notion of a strongly Jacobson ring. 
In classical mathematics, a ring $A$ is strongly Jacobson if and only if $A$ is $\alpha$-Jacobson for some ordinal $\alpha$.
In a foundation of predicative mathematics with generalized inductive definitions, the definition of a strongly Jacobson ring (definition \ref{def-stjac}) is more acceptable than that of a Jacobson ring (definition \ref{def-jac}).
The general Nullstellensatz for strongly Jacobson rings (theorem \ref{qgn-stjac}) also holds.
\section{\texorpdfstring{$\alpha$}{α}-Jacobson rings}
We first recall the definition of a Jacobson ring.
\begin{Def}
    Let ${\gen U}_A$ denote the ideal of a ring $A$ generated by a subset $U\subseteq A$. We define two ideals $\Nil_A U$, $\Jac_A U$ of $A$ as follows:
    \begin{gather*}
    \Nil_A U:=\{x\in A:\exists n\ge0.\ x^n\in\gen U\},\\\Jac_AU:=\{x\in A:\forall a\in A.\ 1\in\gen{U,1-ax}\}.
    \end{gather*}
    Here $\gen{U,1-ax}$ means $\gen{U\cup\{1-ax\}}$. When the context is clear, we write $\gen U$, $\Nil U$, and $\Jac U$ for ${\gen U}_A$, $\Nil_AU$, and $\Jac_AU$, respectively.
    Note that $\Nil U=\Nil{\gen{U}}$ and $\Jac U=\Jac{\gen{U}}$ hold for all $U\subseteq A$.
\end{Def}
\begin{Def}[{\cite[Section 2.4.1]{We18}}]\label{def-jac}
    We call a ring $A$ \emph{Jacobson} if every subset $U$ of $A$ satisfies $\Jac U\subseteq\Nil U$.
\end{Def}
\begin{Exp}
    Every $0$-dimensional ring is Jacobson \cite[Lemma IX-1.2]{LQ15}. The ring $\Z$ is Jacobson \cite[Proposition 2.2]{Ku24}.
    Every algebra that is integral over a Jacobson ring is Jacobson \cite[Corollary 3.8]{Ku24}.
    Every finitely generated algebra over a Jacobson ring is Jacobson \cite[Corollary 3.10]{Ku24}.
\end{Exp}
\begin{Def}
    Let $A$ be a ring, $x,x'\in A$, and $\alpha$ be an ordinal.
    \begin{enumerate}
        \item The triple $(A,x,x')$ is called \emph{$0$-Jacobson} if $x'\in\Nil_A0$.
        \item Let $\alpha>0$. The triple $(A,x,x')$ is called $\alpha$-Jacobson if there exist $n\ge0$ and $a_{1},\ldots,a_{n}\in A$ such that for any $b_{1},\ldots,b_{n}\in A$, there exists $\beta<\alpha$ such that $(A/\gen{1-b_{i}(1-a_{i}x):i\in\{1,\ldots,n\}},x,x')$ is $\beta$-Jacobson.
    \end{enumerate}
    A ring $A$ is called $\alpha$-Jacobson if $(A,x,x)$ is $\alpha$-Jacobson for all $x\in A$.
\end{Def}
\begin{Rem}[ZF]
    We have introduced the notion of an $\alpha$-Jacobson ring to measure the complexity of an elementary proof that a ring is Jacobson.
    In classical mathematics, a triple $(A,x,x')$ is $\alpha$-Jacobson if and only if {\prover} has a winning strategy for the game $J_\alpha(A,x,x')$ defined as follows:
    \begin{enumerate}
        \item Let $\alpha$ be an ordinal, $A$ be a ring, and $x\in A$. The game $J_\alpha(A,x,x')$ is played by two players called {\prover} and {\delayer}.
        \item A possible position of the game is a pair $(\tau,U)$ of an ordinal $\tau\le\alpha$ and a finite subset $U$ of $A$.
        \item The initial position of the game is $(\alpha,\emptyset)$.
        \item Let $(\tau,U)$ be the current position.
        \begin{itemize}
            \item If $\tau>0$, then {\prover} declares a natural number $n\in\N$ and $n$ elements $a_{1},\ldots,a_{n}\in A$. Then {\delayer} declares $n$ elements $b_{1},\ldots,b_{n}\in A$. Then {\prover} declares an ordinal $\tau'<\tau$. The next position is $(\tau',U')$, where $U':=U\cup\{1-b_{i}(1-a_{i}x):i\in\{1,\ldots,n\}\}$.
            \item If $\tau=0$, then the game ends. {\prover} wins if $x'\in\Nil U$. {\delayer} wins if $x'\notin\Nil U$.
        \end{itemize}
    \end{enumerate}
    In the game $J_\alpha(A,x,x')$, {\prover} is trying to give an elementary proof that $(x\in\Jac U)\to(x'\in\Nil U)$ holds for all subsets $U\subseteq A$.
    The idea of using the Prover-Delayer game comes from the theory of proof complexity, which the author has learned from Ken \cite{Ke24}.
    Another game related to the Noetherianity of rings is considered in \cite[Section 5]{Br18}.
\end{Rem}
\begin{Rem}
    In constructive mathematics, we have to be careful what the term \emph{ordinal} means. See \cite{CLN22,KNFX23} for the constructive theory of ordinals. We only need finite ordinals to treat the examples in this paper.
\end{Rem}
\begin{Rem}\label{a-jac-rem}
Let $\alpha,\alpha'$ be ordinals, $A$ be a ring, and $x,x',y,z\in A$.
    \begin{enumerate}
        \item If $\alpha\le\alpha'$ and $(A,x,x')$ is $\alpha$-Jacobson, then $(A,x,x')$ is $\alpha'$-Jacobson.
        \item If $(A,x,x')$ is $\alpha$-Jacobson and $I$ is an ideal of $A$, then $(A/I,x,x')$ is $\alpha$-Jacobson.
        \item If $(A,xy,x')$ is $\alpha$-Jacobson, then $(A,x,x'z)$ is $\alpha$-Jacobson.
        \item The ring $A$ is $0$-Jacobson if and only if $1=_A0$.
    \end{enumerate}
\end{Rem}
We next prove that every $\alpha$-Jacobson ring is Jacobson.
\begin{Th}
    Let $\alpha$ be an ordinal, $A$ be a ring, $x,x'\in A$, and $U\subseteq A$. If $(A,x,x')$ is $\alpha$-Jacobson and $x\in\Jac U$, then $x'\in\Nil U$.
\end{Th}
\begin{Proof}
    We prove this by induction on $\alpha$.
    \begin{enumerate}
        \item If $\alpha=0$, then $x'\in\Nil0\subseteq\Nil U$ holds.
        \item Let $\alpha>0$.
        There exist $n\ge0$ and $a_{1},\ldots,a_{n}\in A$ such that for all $b_{1},\ldots,b_{n}\in A$, there exists $\beta<\alpha$ such that $(A/\gen{1-b_{i}(1-a_{i}x):i\in\{1,\ldots,n\}},x,x')$ is $\beta$-Jacobson.
        Since $x\in\Jac U$, there exist $b_{1},\ldots,b_{n}\in A$ such that $1-b_{i}(1-a_{i}x)\in\gen{U}$ for all $i\in\{1,\ldots,n\}$.
        Hence there exists $\beta<\alpha$ such that $(A/\gen{U},x,x')$ is $\beta$-Jacobson.
        Hence $x'\in\Nil_{A/\gen{U}}0$ by the inductive hypothesis.
        Hence $x'\in\Nil_AU$.\qedhere
    \end{enumerate}
\end{Proof}
\begin{Cor}\label{a-jac-jac}
    Let $\alpha$ be an ordinal. Then every $\alpha$-Jacobson ring is Jacobson.
\end{Cor}
We present some examples of $\alpha$-Jacobson rings.
\begin{Pro}[{a quantitative version of \cite[Example 2.6]{Ku24}}]\label{0-dim-jac}
    Let $A$ be a ring, $x,a\in A$, and $e\ge0$. If $x^e(1-ax)=0$, then $(A,x,x)$ is $1$-Jacobson.
    In particular, every $0$-dimensional ring $A$ is $1$-Jacobson.
\end{Pro}
\begin{Proof}
    Let $b\in A$. Then $x^e=x^e(1-b(1-ax))\in\gen{1-b(1-ax)}$.
\end{Proof}
Let $\Reg A$ denote the set of regular elements of a ring $A$.
\begin{Pro}[{a quantitative version of \cite[Corollary 2.8]{Ku24}}]\label{1-dim-jac}
    Every ring $A$ satisfying the following conditions is $2$-Jacobson:
    \begin{enumerate}
        \item The ring $A$ is integral (i.e., every element is null or regular).
        \item The Krull dimension of $A$ is less than $2$.
        \item For any $x\in\Reg A$, there exists $a\in A$ such that $(1-ax\in A^\times)\to(1=_A0)$.
    \end{enumerate}
    In particular, $\Z$ is $2$-Jacobson, and $K[X]$ is $2$-Jacobson for every discrete field $K$. Note that a ring $K$ is called a discrete field if every element of $K$ is null or invertible.
\end{Pro}
\begin{Proof}
    It suffices to prove that $(A,x,x)$ is $2$-Jacobson for all $x\in\Reg A$.
    Let $x\in\Reg A$. Then there exists $a\in A$ such that $(1-ax\in A^\times)\to(1=_A0)$.
    Let $b\in A$.
    \begin{enumerate}
        \item If $1-b(1-ax)=0$, then $1-ax\in A^\times$. Hence $1=_A0$, and $(A,x,x)$ is $0$-Jacobson.
        \item If $1-b(1-ax)\in\Reg A$, then $A/\gen{1-b(1-ax)}$ is $0$-dimensional. Hence $(A/\gen{1-b(1-ax)},x)$ is $1$-Jacobson by proposition \ref{0-dim-jac}.
        Hence $(A,x,x)$ is $2$-Jacobson.\qedhere
    \end{enumerate}
\end{Proof}
We use the following lemmas to prove that $\Z$ and $A[X]$ are not $1$-Jacobson, where $A$ is a non-trivial ring:
\begin{Lem}[{\cite[Lemma 3.4]{Ku24}}]\label{radical-properties}
    Let $A$ be a ring, $U$ be a subset of $A$, and $x,y\in A$.
    If $xy\in\Nil U$ and $x\in\Nil{(U\cup\{y\})}$, then $x\in\Nil U$.
\end{Lem}
\begin{Lem}[{\cite[Lemma II-2.6]{LQ15}}]\label{invertible-polynomial}
    For all rings $A$, we have $A[X]^\times\subseteq A^\times+(\Nil_A0)[X]$.
\end{Lem}
\begin{Pro}\label{Z-lower-bound}
    If $N$ is an integer such that $|N|\ge2$, the triple $(\Z,N,N)$ is not $1$-Jacobson. In particular, $\Z$ is not $1$-Jacobson.
\end{Pro}
\begin{Proof}
    Assume that $(\Z,N,N)$ is $1$-Jacobson. Then there exist $a_{1},\ldots,a_{n}\in\Z$ such that for all $b_{1},\ldots,b_{n}\in\Z$, we have $N\in\Nil\gen{1-b_{1}(1-a_{1}N),\ldots,1-b_{n}(1-a_{n}N)}$. Let $c:=1+|N(1-a_{1}N)\cdots(1-a_{n}N)|$. Then there exist $b_{1},\ldots,b_{n}\in\Z$ such that $1-b_{i}(1-a_{i}N)=c$ for all $i\in\{1,\ldots,n\}$.
    Hence $N\in\Nil c$. Hence $1\in\Nil c$ by $1\in\gen{c,N}$ and lemma \ref{radical-properties}. This contradicts the fact that $c\ge2$.
\end{Proof}
\begin{Pro}\label{poly-lower-bound}
    Let $A$ be a ring. If $(A[X],X,X)$ is $1$-Jacobson, then $1=_A0$.
\end{Pro}
\begin{Proof}
    Since $(A[X],X,X)$ is $1$-Jacobson, there exist $f_1,\ldots,f_n\in A[X]$ such that for all $g_1,\ldots,g_n\in A[X]$, we have $X\in\Nil\gen{1-g_1(1-f_1X),\ldots,1-g_n(1-f_nX)}$. Let $h:=1-X(1-f_1X)\cdots(1-f_nX)$. Then there exist $g_1,\ldots,g_n\in\Z$ such that $1-g_i(1-f_iX)=h$ for all $i\in\{1,\ldots,n\}$.
    Hence $X\in\Nil h$.
    Hence $1\in\Nil h$ by $1\in\gen{X,h}$ and lemma \ref{radical-properties}.
    Hence $h\in A[X]^\times\subseteq A^\times+(\Nil_A0)[X]$ by lemma \ref{invertible-polynomial}.
    Hence $1=_A0$.
\end{Proof}
We next provide a Jacobson ring that is not $\alpha$-Jacobson for any $\alpha$.
For a ring $A$, let $A[X_k:k\in\N]$ denote the polynomial ring in countably infinitely many variables over $A$.
\begin{Lem}\label{lemma-polynomial-not-alpha-jacobson}
    Let $K$ be a ring, $A:=K[X_k:k\in\N]$, $\alpha$ be an ordinal, $n\ge0$, and $f_1,\ldots,f_n\in A$.
    Let $l:\{1,\ldots,n\}\to\N_{\ge1}$ be a function.
    Assume $f_i\in K[X_0,\ldots,X_{l(i)-1}]$ and $\forall j<i.\ l(j)<l(i)$ for all $i\in\{1,\ldots,n\}$.
    Let $U:=\{1-X_{l(i)}(1-f_iX_0):i\in\{1,\ldots,n\}\}$.
    If $(A/\gen{U},X_0,X_0)$ is $\alpha$-Jacobson, then $1=_K0$.
\end{Lem}
\begin{Proof}
    We prove this by induction on $\alpha$.
    \begin{enumerate}
        \item If $(A/\gen{U},X_0,X_0)$ is $0$-Jacobson, then $X_0\in\Nil_{A/\gen{U}}0$. Let $S:=1+XK[X]\subseteq K[X]$ and $B:=S^{-1}(K[X])$.
        We define a $K$-homomorphism $\varphi:A/\gen{U}\to B$ by the following equation:
        \[
        \varphi(X_k):=\begin{cases}
            X&\text{if $\forall i\in\{1,\ldots,n\}.\ k\ne l(i)$,}\\
            (1-\varphi(f_i)X)^{-1}&\text{if $k=l(i)$.}
        \end{cases}
        \]
        Then $X=\varphi(X_0)\in\Nil_B0$.
        Since $K[X]\subseteq B$, we have $1=_K0$.
        \item Let $\alpha>0$. Since $(A/\gen{U},X_0,X_0)$ is $\alpha$-Jacobson, there exist $m>n$ and $f_{n+1},\ldots,f_m\in A$ such that for all $g_{n+1},\ldots,g_m\in A$, there exists $\beta<\alpha$ such that \[(A/\gen{U\cup\{1-g_i(1-f_iX_0):i\in\{n+1,\ldots,m\}\}},X_0,X_0)\] is $\beta$-Jacobson.
        There exist $l(n+1),\ldots,l(m)\in\N$ such that $f_i\in K[X_0,\ldots,X_{l(i)-1}]$ and $\forall j<i.\ l(j)<l(i)$ hold for every $i\in\{1,\ldots,m\}$.
        Hence there exists $\beta<\alpha$ such that
        \[(A/\gen{1-X_{l(i)}(1-f_{i}X_0):i\in\{1,\ldots,m\}},X_0,X_0)\] is $\beta$-Jacobson.
        Hence $1=_K0$ by the inductive hypothesis.\qedhere
    \end{enumerate}
\end{Proof}
\begin{Cor}\label{polynomial-not-alpha-jacobson}
    Let $K$ be a ring. If $K[X_k:k\in\N]$ is $\alpha$-Jacobson for some $\alpha$, then $1=_K0$.
\end{Cor}
\begin{Rem}[ZFC]\label{not-a-jac}
    Krull \cite[Satz 4]{Kr51} has proved that if $K$ is an uncountable field, then $K[X_k:k\in\N]$ is Jacobson. This also follows from Amitsur's theorem \cite[Corollary 3]{Am56}.
    Hence $\R[X_k:k\in\N]$ is Jacobson, but it is not $\alpha$-Jacobson for any $\alpha$ by corollary \ref{polynomial-not-alpha-jacobson}.
\end{Rem}
\section{The quantitative general Nullstellensatz}
In this section, we prove the quantitative general Nullstellensatz: if $A$ is an $\alpha$-Jacobson ring, $A[X]$ is $(\alpha+1)$-Jacobson.
The proof is similar to the constructive proof of the general Nullstellensatz, but we need more exact lemmas than those in \cite{Ku24}.
\begin{Lem}\label{cut}
    Let $A$ be a ring and $x,y,z\in A$.
    If $(A,y,xz)$ and $(A/\gen{x},y,z)$ are $\alpha$-Jacobson, then $(A,y,z)$ is $\alpha$-Jacobson.
\end{Lem}
\begin{Proof}
    We prove this by induction on $\alpha$.
    \begin{enumerate}
        \item If $\alpha=0$, then $xz\in\Nil_A0$ and $z\in\Nil_Ax$. Hence $z\in\Nil_A0$ by lemma \ref{radical-properties}.
        \item Let $\alpha>0$. Then there exist $m,n\ge0$, $a_1,\ldots,a_m\in A$, and $a_1',\ldots,a_n'\in A$ such that for all $b_1,\ldots,b_m\in A$ and $b_1',\ldots,b_n'\in A$, there exist ordinals $\beta,\beta'$ such that
        \begin{itemize}
            \item $(A/\gen{1-b_i(1-a_iy):i\in\{1,\ldots,m\}},y,xz)$ is $\beta$-Jacobson, and
            \item $(A/(\gen{x}+\gen{1-b_i'(1-a_i'y):i\in\{1,\ldots,n\}}),y,z)$ is $\beta'$-Jacobson.
        \end{itemize}
        Note that the above two items together imply that
        \begin{itemize}
            \item $(A/(\gen{1-b_i(1-a_iy):i\in\{1,\ldots,m\}}+\gen{1-b_i'(1-a_i'y):i\in\{1,\ldots,n\}}),y,z)$ is $\gamma$-Jacobson, where $\gamma:=\max(\beta,\beta')$.
        \end{itemize}
        by the inductive hypothesis.
        Hence $(A,y,z)$ is $\alpha$-Jacobson.\qedhere
    \end{enumerate}
\end{Proof}
\begin{Lem}\label{integral-entail}
    Let $A$ be a ring, $B$ be an $A$-algebra which is integral over $A$, and $x\in A$. Then
    \[\forall b\in B.\ \exists a\in A.\ 1-ax\in\gen{1-bx}_B.\]
\end{Lem}
\begin{Proof}
    Let $b\in B$.
    Since $C:=B/\gen{1-bx}$ is integral over $A$ and $x\in C^\times$, there exists $a\in A$ such that $ax=_C1$ by \cite[Theorem IX-1.7]{LQ15}. Hence we have $1-ax\in\gen{1-bx}_B$.
\end{Proof}
\begin{Lem}\label{loc-key}
    Let $A$ be a ring, $a,a_1,a_2'\in A$, and $e\ge0$.
    Then there exists $a_2\in A$ such that
    \[1-a_2(1-a_1a)\in\gen{a^e-a_2'(1-a_1a)}_A.\]
\end{Lem}
\begin{Proof}
    Let $a_2:=(1+\cdots+(a_1a)^{e-1})+a_1^ea_2'$.
    Then,
    \begin{align*}
        1-a_2(1-a_1a)&=1-((1+\cdots+(a_1a)^{e-1})+a_1^ea_2')(1-a_1a)\\
        &=(a_1a)^e-a_1^ea_2'(1-a_1a)\\
        &=a_1^e(a^e-a_2'(1-a_1a)).
    \end{align*}
    Hence $1-a_2(1-a_1a)\in\gen{a^e-a_2'(1-a_1a)}_A$.
\end{Proof}
For a ring $A$ and an element $a\in A$, let $A_a$ denote the ring $A[1/a]$.
\begin{Lem}\label{key-elementary}
    Let $A$ be a ring, $a,a_0\in A$, $B$ be an $A$-algebra such that $B_a$ is integral over $A_a$. Then
        \[\forall a_1\in A.\ \forall b_2\in B.\ \exists a_2\in A.\ 1-a_2(1-a_1aa_0)\in\gen{1-b_2(1-a_1aa_0)}_B.\]
\end{Lem}
\begin{Proof}
    Let $a_1\in A$ and $b_2\in B$. By lemma \ref{integral-entail}, there exists $a_2'\in A_a$ such that
    \[
    1-a_2'(1-a_1aa_0)\in\gen{1-b_2(1-a_1aa_0)}_{B_a}.
    \]
    There exists $a_2''\in A$ and $e\ge0$ such that $a_2'=a_2''/a^e$.
    Then there exists $e'\ge0$ such that
    \[
    a^{e'}(a^e-a_2''(1-a_1aa_0))\in\gen{1-b_2(1-a_1aa_0)}_B.
    \]
    By lemma \ref{loc-key}, there exists $a_2\in A$ such that
    \[
    1-a_2(1-a_1aa_0)\in\gen{a^{e+e'}-a^{e'}a_2''(1-a_1aa_0)}_A.
    \]
    Hence
    \[
    1-a_2(1-a_1aa_0)\in\gen{1-b_2(1-a_1aa_0)}_B
    \]
    holds.
\end{Proof}
The following lemma follows from the definition:
\begin{Lem}\label{cumbersome}
    Let $(A,x,x')$ be an $\alpha$-Jacobson triple, $B$ be an $A$-algebra, and $y\in B$. If 
    \[\forall a_1\in A.\ \exists b_1\in B.\ \forall b_2\in B.\ \exists a_2\in A.\ 1-a_2(1-a_1x)\in\gen{1-b_2(1-b_1y)}_B,\]
    then there exist $n\ge0$, $a_{1,1},\ldots,a_{1,n}\in A$, and $b_{1,1},\ldots,b_{1,n}\in B$ such that for all $b_{2,1},\ldots,b_{2,n}\in B$, there exists $a_{2,1},\ldots,a_{2,n}\in A$ such that
    \begin{itemize}
        \item there exists $\beta<\alpha$ such that $(A/\gen{1-a_{2,i}(1-a_{1,i}x):i\in\{1,\ldots,n\}},x,x')$ is $\beta$-Jacobson, and
        \item $1-a_{2,i}(1-a_{1,i}x)\in\gen{1-b_{2,i}(1-b_{1,i}y)}_B$ holds for all $i\in\{1,\ldots,n\}$.
    \end{itemize}
\end{Lem}
\begin{Lem}\label{loc-integral-jac}
    Let $\alpha$ be an ordinal, $A$ be a ring, $a\in A$, $B$ be an $A$-algebra such that $B_a$ is integral over $A_a$, and $y\in B$. If there exist $d\ge0$, $c_0,\ldots,c_{d-1}\in A$, and $l\ge0$ such that $a^ly^d=_Bc_{d-1}y^{d-1}+\cdots+c_0$ and $(A,ac_k,ac_k)$ is $\alpha$-Jacobson for all $k\in\{0,\ldots,d-1\}$, then $(B,y,ay)$ is $\alpha$-Jacobson.
\end{Lem}
\begin{Proof}
    We prove this by induction on $\alpha$.
    \begin{enumerate}
        \item If $\alpha=0$, then $ac_0,\ldots,ac_{d-1}\in\Nil_A0$. Hence $ay\in\Nil_B0$, and $(B,y,ay)$ is $0$-Jacobson.
        \item Let $\alpha>0$.
        Let $f_k:=a^ly^k-(c_{d-1}y^{k-1}+\cdots+c_{d-k})$ and $B_k:=B/\gen{f_k,\ldots,f_{d-1}}$ for $k\in\{0,\ldots,d-1\}$.
        We prove that $(B_k,y,ay)$ is $\alpha$-Jacobson for all $k$ by induction.
        \begin{enumerate}
            \item Since $a^l=f_0=_{B_0}0$, the triple $(B_0,y,ay)$ is $0$-Jacobson.
            \item Let $k\ge1$.
        \begin{itemize}
            \item 
            By lemma \ref{key-elementary} and lemma \ref{cumbersome}, there exist $n\ge0$, $a_{1,1},\ldots,a_{1,n}\in A$, and $b_{1,1},\ldots,b_{1,n}\in B_k$ such that for all $b_{2,1},\ldots,b_{2,n}\in B_k$, there exist $a_{2,1},\ldots,a_{2,n}\in A$ such that
    \begin{itemize}
        \item there exists $\beta<\alpha$ such that \[(A/\gen{1-a_{2,i}(1-a_{1,i}ac_{d-k}):i\in\{1,\ldots,n\}},ac_{d-k},ac_{d-k})\] is $\beta$-Jacobson, and
        \item $1-a_{2,i}(1-a_{1,i}ac_{d-k})\in\gen{1-b_{2,i}(1-b_{1,i}ac_{d-k})}_B$ holds for all $i\in\{1,\ldots,n\}$.
    \end{itemize}
    The above two items together imply that
    \begin{itemize}
        \item there exists $\beta<\alpha$ such that \[(B_k/\gen{1-b_{2,i}(1-(b_{1,i}a)c_{d-k}):i\in\{1,\ldots,n\}},c_{d-k},ac_{d-k})\] is $\beta$-Jacobson
    \end{itemize}
    by the inductive hypothesis.
    Hence $(B_k,c_{d-k},ac_{d-k})$ is $\alpha$-Jacobson.
    Since $c_{d-k}=f_{k-1}y-f_k=_{B_k}f_{k-1}y$, the triple $(B_k,y,af_{k-1}y)$ is $\alpha$-Jacobson by remark \ref{a-jac-rem}-3.
    \item The triple $(B_{k-1},y,ay)$ is $\alpha$-Jacobson by the inductive hypothesis.
    \end{itemize}
    Hence $(B_k,y,ay)$ is $\alpha$-Jacobson by lemma \ref{cut}.\qedhere
    \end{enumerate}
    \end{enumerate}
\end{Proof}
The above lemma implies the following theorem:
\begin{Th}[{a quantitative version of \cite[Lemma 3.6]{Ku24}}]
    Let $\alpha$ be an ordinal, $A$ be an $\alpha$-Jacobson ring, $a\in A$, $B$ be an $A$-algebra such that $B_a$ is integral over $A_a$, and $y\in B$. Then $(B,y,ay)$ is $\alpha$-Jacobson.
\end{Th}
\begin{Cor}[a quantitative version of {\cite[Corollary 3.7]{Ku24}}]
    Let $\alpha$ be an ordinal, $A$ be an $\alpha$-Jacobson ring, and $a\in A$. Then $A_a$ is $\alpha$-Jacobson.
\end{Cor}
\begin{Cor}[a quantitative version of {\cite[Corollary 3.8]{Ku24}}]
    Let $\alpha$ be an ordinal, $A$ be an $\alpha$-Jacobson ring. Then every $A$-algebra $B$ that is integral over $A$ is $\alpha$-Jacobson.
\end{Cor}
Using lemma \ref{loc-integral-jac}, we prove the quantitative general Nullstellensatz, which is a quantitative version of {\cite[Theorem 3.9]{Ku24}}.
We first recall the following lemma:
\begin{Lem}[{\cite[Corollary VI-1.3]{MRR88}}]\label{integrality-monogenic-algebra}
    Let $B$ be an $A$-algebra. If $b\in B$ is integral over $A$, then the algebra $A[b]\subseteq B$ is integral over $A$.
\end{Lem}
\begin{Th}
    Let $\alpha$ be an ordinal, $A$ be an $\alpha$-Jacobson ring, and $f,g\in A[X]$. Then
    \[(A[X]/\gen{1-g(1-Xf)},f,f)\] is $\alpha$-Jacobson.
\end{Th}
\begin{Proof}
    Let $h:=1-g(1-Xf)$. There exist $d\ge0$ and $a_0,\ldots,a_d\in A$ such that $h=a_dX^d+\cdots+a_0$.
    Let $C_k:=A[X]/\gen{h,a_{k+1},\ldots,a_d}$ for $k\in\{-1,\ldots,d\}$.
    We prove that $(C_k,f,f)$ is $\alpha$-Jacobson by induction on $k$.
    \begin{enumerate}
        \item Since $h=_{C_{-1}}0$, we have $1-Xf\in C_{-1}^\times$.
        Since $C_{-1}=(A/\gen{a_0,\ldots,a_d})[X]$, we have $f\in\Nil_{C_{-1}}0$ by lemma \ref{invertible-polynomial}.
        Hence $(C_{-1},f,f)$ is $0$-Jacobson.
        \item Let $k\ge0$.
        \begin{itemize}
            \item The $A_{a_k}$-algebra $(C_k)_{a_k}$ is integral over $A_{a_k}$ by lemma \ref{integrality-monogenic-algebra}. Hence the triple $(C_k,f,a_kf)$ is $\alpha$-Jacobson by lemma \ref{loc-integral-jac}.
            \item The triple $(C_{k-1},f,f)$ is $\alpha$-Jacobson by the inductive hypothesis.
        \end{itemize}
        Hence $(C_k,f,f)$ is $\alpha$-Jacobson by lemma \ref{cut}.\qedhere
    \end{enumerate}
\end{Proof}
\begin{Cor}[The quantitative general Nullstellensatz]\label{qgn}
    If $A$ is an $\alpha$-Jacobson ring, then $A[X]$ is $(\alpha+1)$-Jacobson.
\end{Cor}
In \cite{Ku24}, we have proved that $K[X_1,\ldots,X_n]$ is Jacobson for any discrete field $K$ and that $\mathbb{Z}[X_1,\ldots,X_n]$ is Jacobson. These results provide a solution to the first two problems of Lombardi's list \cite{Lo23}.
By proposition \ref{0-dim-jac}, proposition \ref{1-dim-jac}, and corollary \ref{qgn}, we have the following stronger results: $K[X_1,\ldots,X_n]$ is $(1+n)$-Jacobson for any discrete field $K$, and $\mathbb{Z}[X_1,\ldots,X_n]$ is $(2+n)$-Jacobson.
\section{Strongly Jacobson rings}
In this section, we introduce the notion of strongly Jacobson rings.
\begin{Def}\label{def-stjac}
    Let $A$ be a ring and $x,x'\in A$.
    The proposition ``The triple $(A,x,x')$ is \emph{strongly Jacobson}'' is inductively generated by the following constructors:
    \begin{enumerate}
        \item If $x'\in\Nil_A0$, then $(A,x,x')$ is strongly Jacobson.
        \item If there exist $n\ge0$ and $a_{1,1},\ldots,a_{1,n}\in A$ such that for any $a_{2,1},\ldots,a_{2,n}\in A$, the triple \[(A/\gen{1-a_{2,i}(1-a_{1,i}x):i\in\{1,\ldots,n\}},x,x')\] is strongly Jacobson,
        then $(A,x,x')$ is strongly Jacobson.
    \end{enumerate}
    A ring $A$ is called strongly Jacobson if $(A,x,x)$ is strongly Jacobson for all $x\in A$.
\end{Def}
\begin{Rem}
    Definition \ref{def-stjac} is a generalized inductive definition.
    We have simultaneously defined the proposition ``The triple $(A/I,x,x')$ is strongly Jacobson'' for all finitely generated ideals $I$ of $A$.
    Generalized inductive definitions are accepted in some foundations of predicative mathematics. They are also used to define a well-founded relation \cite[Section 10.3]{UFP13} and a Noetherian ring \cite{CP99,CL01,BSB23}.
\end{Rem}
The following theorem easily follows from the definition.
\begin{Th}
    Let $A$ be an $\alpha$-Jacobson ring for some $\alpha$. Then $A$ is strongly Jacobson.
\end{Th}
\begin{Th}[ZF]
    Let $(A,x,x')$ be a strongly Jacobson triple. Then there exists an ordinal $\alpha$ such that $(A,x,x')$ is $\alpha$-Jacobson.
\end{Th}
\begin{Proof}
    We prove this by induction on the proof of the proposition ``The triple $(A,x,x')$ is strongly Jacobson.''
    \begin{enumerate}
        \item If $x'\in\Nil_A0$, then $(A,x,x')$ is $0$-Jacobson.
        \item Assume that there exist $n\ge0$ and $a_{1,1},\ldots,a_{1,n}\in A$ such that for any $a_{2,1},\ldots,a_{2,n}\in A$, the triple \[(A/\gen{1-a_{2,i}(1-a_{1,i}x):i\in\{1,\ldots,n\}},x,x')\] is strongly Jacobson. By the inductive hypothesis, for any $a_{2,1},\ldots,a_{2,n}\in A$, there exists the minimum ordinal $\beta_{a_{2,1},\ldots,a_{2,n}}$ such that $(A/\gen{1-a_{2,i}(1-a_{1,i}x):i\in\{1,\ldots,n\}},x,x')$ is $(\beta_{a_{2,1},\ldots,a_{2,n}})$-Jacobson. Then $A$ is $(\sup_{a_{2,1},\ldots,a_{2,n}}\beta_{a_{2,1},\ldots,a_{2,n}}+1)$-Jacobson.\qedhere
    \end{enumerate}
\end{Proof}
We can prove the following theorem by an argument similar to corollary \ref{a-jac-jac}.
\begin{Th}
     Every strongly Jacobson ring is Jacobson.
\end{Th}
The converse does not hold by remark \ref{not-a-jac}.
We can prove the following theorem by an argument similar to corollary \ref{qgn}.
\begin{Th}[The general Nullstellensatz for strongly Jacobson rings]\label{qgn-stjac}
    If $A$ is a strongly Jacobson ring, then $A[X]$ is strongly Jacobson.
\end{Th}
\section{Future work}
We present three problems on $\alpha$-Jacobson rings and strongly Jacobson rings.
\begin{Prob}
    Let $K$ be a field.
    What is the minimum ordinal $\alpha$ such that $K[X_1,\ldots,X_n]$ is $\alpha$-Jacobson?
\end{Prob}
Corollary \ref{qgn} implies that $K[X_1,\ldots,X_n]$ is $(1+n)$-Jacobson, and proposition \ref{poly-lower-bound} implies that if $K[X]$ is $0$-Jacobson, then $1=_K0$.
\begin{Prob}
    What is the minimum ordinal $\alpha$ such that $\Z[X_1,\ldots,X_n]$ is $\alpha$-Jacobson?
\end{Prob}
Corollary \ref{qgn} implies that $\Z[X_1,\ldots,X_n]$ is $(2+n)$-Jacobson, and proposition \ref{Z-lower-bound} implies that $\Z$ is not $1$-Jacobson.
\begin{Prob}
    Is there a Noetherian Jacobson ring that is not strongly Jacobson?
\end{Prob}
The ring $\R[X_k:k\in\N]$ is a non-Noetherian Jacobson ring that is not strongly Jacobson by remark \ref{not-a-jac}.
\section*{Acknowledgments}
The author would like to express his deepest gratitude to his supervisor, Ryu Hasegawa, for his support.
The author would like to thank Eitetsu Ken for the inspiring PhD thesis presentation, which led the author to the definition of an $\alpha$-Jacobson ring.
The author would also like to thank Yuto Ikeda for the interesting discussion.

This research was supported by Forefront Physics and Mathematics Program to Drive Transformation (FoPM), a World-leading Innovative Graduate Study (WINGS) Program, the University of Tokyo.


\begin{thebibliography}{99}
    \bibitem[Ami56]{Am56} A.~S.~Amitsur. Algebras over infinite fields, \emph{Proc.~Amer.~Math.~Soc.} 7:35--48 (1956).
    \bibitem[Bra18]{Br18} M.~Brandenburg. Algebraic games---playing with groups and rings. \emph{Int.~J.~Game Theory} 47:417--450 (2018).
    \bibitem[BSB23]{BSB23} G.~Buriola, P.~Schuster, I.~Blechschmidt. A Constructive Picture of Noetherian Conditions and Well Quasi-orders. In: G.~Della Vedova, B.~Dundua, S.~Lempp, F.~Manea, eds., \emph{Unity of Logic and Computation}, 50--62. CiE 2023. Lecture Notes in Comput. Sci., Vol.~13967. Springer, Cham (2023).
    \bibitem[CL01]{CL01} T.~Coquand, H.~Lombardi. Krull's Principal Ideal Theorem. Tech. Rep. 30, Institut Mittag-Leffler (2001).
    \bibitem[CLN22]{CLN22} T.~Coquand, H.~Lombardi, S.~Neuwirth, Constructive theory of ordinals. In: M.~Benini, O.~Beyersdorff, M.~Rathjen, P.~Schuster, eds., \emph{Mathematics for Computation (M4C)}. World Scientific, pp.~287--318 (2022).
    \bibitem[CP99]{CP99} T.~Coquand, H.~Persson. Gröbner bases in type theory. In: T.~Altenkirch, B.~Reus, W.~Naraschewski, eds., \emph{TYPES 1998: Types for Proofs and Programs}. Lecture Notes in Comput. Sci., Vol. 1657, pp.~33–46. Springer (1999).
    \bibitem[Eme]{Em} M.~Emerton. \emph{Jacobson rings}. Available at: \url{https://www.math.uchicago.edu/~emerton/pdffiles/jacobson.pdf}.
    \bibitem[Gol51]{Go51} O.~Goldman. Hilbert rings and the Hilbert Nullstellensatz. \emph{Math.~Z.} 54:136–140 (1951).
    \bibitem[Ken24]{Ke24} E.~Ken. \emph{Games with backtracking options corresponding to the ordinal analysis of $PA$}, {\tt arXiv:2406.17315v2 [math.LO]} (2024).
    \bibitem[KNFX23]{KNFX23} N.~Kraus, F.~Nordvall Forsberg, C.~Xu. Type-theoretic approaches to ordinals. \emph{Theoret. Comput. Sci.} 957, Paper No.~113843 (2023).
    \bibitem[Kru51]{Kr51} W.~Krull. Jacobsonsche Ringe, Hilbertscher Nullstellensatz, Dimensionstheorie. \emph{Math.~Z.} 54:354--387 (1951).
    \bibitem[Kur24]{Ku24} R.~Kuroki. \emph{A constructive proof of the general Nullstellensatz for Jacobson rings}, {\tt arXiv:2406.06078v2 [math.AC]} (2024).
    \bibitem[Lom23]{Lo23} H.~Lombardi. \emph{Some classical results in Algebra needing one or several constructive versions}. Available at: \url{https://groups.google.com/g/constructivenews/c/Z6ZEmRdep8o/m/TNVpuihzAAAJ} (2023).
    \bibitem[LQ15]{LQ15} H.~Lombardi, C.~Quitt\'e. \emph{Commutative algebra: constructive methods}. Algebra and Applications, Vol. 20. Springer, Dordrecht (2015).
    \bibitem[MRR88]{MRR88} R.~Mines, F.~Richman, W.~Ruitenburg. \emph{A course in constructive algebra}. Universitext. Springer-Verlag, New York (1988).
    \bibitem[UFP13]{UFP13} The Univalent Foundations Program. \emph{Homotopy Type Theory: Univalent Foundations of Mathematics}. Institute for Advanced Study. Available at: \url{https://homotopytypetheory.org/book} (2013).
    \bibitem[Wes18]{We18} D.~Wessel. \emph{Choice, extension, conservation. From transfinite to finite proof methods in abstract algebra}. PhD thesis. University of Trento, University of Verona (2018).
\end{thebibliography}
\end{document}